\documentclass[a4paper]{article}
\usepackage[AMS,theorem]{cjssac}
\usepackage{geometry}

\usepackage{algorithm}%
\usepackage{algorithmicx}%
\usepackage{algpseudocode}%
\usepackage{enumerate}
\usepackage{cleveref}
\crefname{Th}{Theorem}{Theorems}
\usepackage{bm}

\DeclareMathOperator{\Bez}{Bez}

\title{The GPGCD Algorithm with the B\'ezout Matrix\\ for Multiple Univariate Polynomials}

\author{Boming Chi
\mail{hakumei-t@math.tsukuba.ac.jp}
\affil{Graduate School of Pure and Applied Sciences, University of Tsukuba}
\and
Akira Terui
\mail{terui@math.tsukuba.ac.jp}
\affil{Faculty of Pure and Applied Sciences, University of Tsukuba}
}
\begin{document}
\begin{abstract}
    We propose a modification of the GPGCD algorithm, which has been presented in our previous research, for calculating approximate greatest common divisor (GCD) of more than 2 univariate polynomials with real coefficients and a given degree. 
    In transferring the approximate GCD problem to a constrained minimization problem, different from the original GPGCD algorithm for multiple polynomials which uses the Sylvester subresultant matrix, the proposed algorithm uses the B\'ezout matrix.
    Experiments show that the proposed algorithm is more efficient than the original GPGCD algorithm for multiple polynomials with maintaining almost the same accuracy for most of the cases.
\end{abstract}


\section{Introduction}\label{sec:Introduction}

With the progress of algebraic computation with polynomials and matrices, we are paying more attention to approximate algebraic algorithms.
Algorithms for calculating approximate GCD, which are approximate algebraic algorithms, 
consider a pair of given polynomials $F$ and $G$ that are relatively prime in general, 
and find $\tilde{F}$ and $\tilde{G}$ which are close to $F$ and $G$, respectively,
in the sense of polynomial norm, and have the greatest common divisor of a certain degree.
These algorithms can be classified into two categories: 1) for a given tolerance (magnitude)
of $\|F-\tilde{F}\|$ and $\|G-\tilde{G}\|$, maximize the degree of approximate GCD, and 2) for a given degree $d$, minimize the magnitude of $\|F-\tilde{F}\|$ and $\|G-\tilde{G}\|$.

In both categories, algorithms based on various methods have been proposed including the Euclidean algorithm (\cite{bec-lab1998b}, \cite{sas-nod89}, \cite{sch1985}), low-rank approximation of the Sylvester matrix or subresultant matrices (\cite{cor-gia-tra-wat1995}, \cite{cor-wat-zhi2004},
\cite{emi-gal-lom1997}, \cite{kal-yan-zhi2006}, \cite{kal-yan-zhi2007}, \cite{sch-spa2016}, \cite{Terui2009}, \cite{TERUI2013127}, \cite{zar-ma-fai2000}, \cite{zhi2003}), Pad\'e approximation (\cite{pan2001b}), and optimizations (\cite{chi-cor-cor1998}, \cite{kar-lak1998}, \cite{10.1007/978-3-030-85165-1_16}, \cite{Sun2007}, \cite{zen2011}).
Among them, the second author of the present paper has proposed the GPGCD algorithm based on low-rank approximation of subresultant
matrices by optimization (\cite{Terui2009}, \cite{Terui2010}, \cite{TERUI2013127}), which belongs to the second category above.
Based on the researches mentioned above, the authors of the present paper have proposed the GPGCD algorithm using the B\'ezout matrix (\cite{10.1007/978-3-030-60026-6_10}), which is the previous research of this paper. 

In this paper, we propose the GPGCD algorithm with the B\'ezout matrix
for multiple polynomials, 
while subresultant matrices have been used in the original GPGCD algorithm 
for multiple polynomials (\cite{Terui2010}). 
We show that the proposed algorithm
is more efficient than the original one 
with maintaining almost the same accuracy for most of the cases.

The paper is organized as follows. 
In \Cref{sec:AppGCDprob}, we give a definition of the approximate GCD problem.
In \Cref{sec:trans}, we give a formulation of the transformation of the approximate GCD problem to the optimization problem using the B\'ezout matrix. 
In \Cref{sec:newton}, we review the modified Newton method used for optimization. 
In \Cref{sec:algorithm}, we illustrate the proposed algorithm and give a running time analysis.
In \Cref{sec:experiments}, the results of experiments are shown. 

\section{Approximate GCD Problem}\label{sec:AppGCDprob}

Let $F_1(x), \dots, F_n(x)$ be univariate polynomials with real coefficients of degree at most $m$:
\begin{equation}\label{eq:defAppGCD}
    F_i(x)=f_{im} x^m+\dots+f_{i0} x^0, \quad i=1, \dots, n, \quad f_{1m} \neq 0.
\end{equation}

In this paper, for a polynomial $F(x)=f_m x^m+\dots+f_0 x^0$, the norm $\|F\|$ denotes the 2-norm defined as $\|F\|_2:=(f_m^2+f_{m-1}^2+\dots+f_0^2)^{\frac12}$.
For a vector $(a_1, \dots, a_n) \in \mathbb{R}^n$, the norm $\|(a_1, \dots, a_n)\|$ denotes the 2-norm defined as $\|(a_1, \dots, a_n)\|_2:=(a_1^2+\dots+a_n^2)^{\frac12}$.

Here, we give a definition of an approximate GCD.

\begin{Def}[Approximate GCD]
  \label{def:AppGCD}
  For polynomials $F_1(x), \dots, F_n(x)$ which are relatively prime in general,
  a positive integer $d$ satisfying that $d<m$, and a positive real number $\epsilon$,  
  if there exist polynomials $\tilde{F}_1(x), \dots, \tilde{F}_n(x)$ such that
  they have a certain GCD as
  \begin{equation}
    \label{eq:approximate-gcd}
    \tilde{F}_i(x)= \bar{F}_i(x)\times \tilde{H}(x), \quad i=1, \dots, n,    
  \end{equation}
  where $\tilde{H}(x)$ is a polynomial of degree $d$, satisfying that 
  $\big{\|}(\|\tilde{F}_1-F_1\|, \dots, \|\tilde{F}_n-F_n\|)\big{\|}<\epsilon$, 
  we call $\tilde{H}(x)$ an approximate GCD of polynomials $F_1(x), \dots, F_n(x)$ with tolerance $\epsilon$.
\end{Def}

Algorithms for calculating approximate GCD can be classified into two categories: 1) for a given tolerance $\epsilon$,
make the degree of approximate GCD as large as possible, and 2) for a given degree $d$, minimize the magnitude $\big{\|}(\|F_1-\tilde{F}_1\|, \dots, \|F_n-\tilde{F}_n\|)\big{\|}$.
In this paper, we focus on the second category of the approximate GCD algorithms, solving the following problem.

\begin{Prob}[Approximate GCD problem]
  \label{prob:AppGCD}
  For given univariate polynomials $F_1(x), \dots, F_n(x)$ as shown in \eqref{eq:defAppGCD} and
  a positive integer $d<m$, 
  find polynomials $\tilde{F}_1(x), \dots, \tilde{F}_n(x)$ and $\tilde{H}(x)$
  satisfying \eqref{eq:approximate-gcd}, with making the perturbation 
  \begin{equation}
      \label{eq:perturbation}
      \varDelta :=\sqrt{\sum_{i=1}^n\|F_i(x)-\tilde{F}_i(x)\|^2}
  \end{equation}
  as small as possible.
\end{Prob}

\section{Transformation of the approximate GCD problem}\label{sec:trans}

For solving the approximate GCD problem, we transfer the approximate GCD problem to a constrained minimization problem with the B\'{e}zout matrix.
\begin{Def}[B\'ezout Matrix \cite{Diaz-Toca2002}]\label{def:Bezout}
    Let $F(x)$ and $G(x)$ be two real polynomials with the degree at most $m$. Then, the matrix  $\Bez(F,G)=(b_{ij})_{i,j=1, \dots, m}$, where
    \[
    \frac{F(x)G(y)-F(y)G(x)}{x-y}=\sum_{i,j=1}^{m}b_{ij}x^{i-1}y^{j-1},
    \]
    is called the B\'ezout matrix associated to $F(x)$ and $G(x)$.
    For polynomials $F_1(x), \dots, F_n(x)$, the matrix
    \[
        \Bez(F_1, \dots, F_n)=
        \begin{pmatrix}
            \Bez(F_1, F_2)\\
            \Bez(F_1, F_3)\\
            \vdots\\
            \Bez(F_1, F_n)\\
        \end{pmatrix}
    \]
    is called the B\'{e}zout matrix associated to $F_1(x), \dots, F_n(x)$.
\end{Def}

For the B\'{e}zout matrix, we have a following theorem.

\begin{Th}[Barnett's theorem \cite{Diaz-Toca2002}]
    \label{Theo:Barnett}
    Let $F_1(x), \dots, F_n(x)$ be real polynomials with the degree at most $m$. Let $d=\deg(\gcd(F_1(x), \dots, F_n(x)))$ and $(\bm{b}_1, \dots, \bm{b}_m)=\Bez(F_1(x), \dots, F_n(x))$.
    Then, the vectors $\bm{b}_{1},\dots,\bm{b}_{m-d}$ are linearly independent, and there exists coefficients $c_{i,j}$ such that
    \begin{equation}
    \bm{b}_i=\sum_{j=1}^{m-d}c_{i,j}\bm{b}_{j},\quad1 \le i \le d,
    \end{equation}
    Furthermore, the monic form of the GCD of $F_1(x), \dots, F_n(x)$ is represented as
    \begin{equation}
    \gcd(F_1(x), \dots, F_n(x))=x^d + c_{d,1}x^{d-1} + \dots + c_{1,1}x^0.
    \end{equation}
\end{Th}

For polynomials $\tilde{F}_1(x), \dots, \tilde{F}_n(x)$ in Problem \ref{prob:AppGCD}, let $\tilde{B}=\Bez(\tilde{F}_1(x), \dots, \tilde{F}_n(x))=(\bm{b}_1, \dots, \bm{b}_m)$.
In the case that polynomials $\tilde{F}_1(x), \dots, \tilde{F}_n(x)$ have an exact GCD of degree $d$, 
Theorem \ref{Theo:Barnett} shows that $\bm{b}_{1}, \dots, \bm{b}_{m-d}$ are linearly independent, 
and $\bm{b}_{m-d+1}$ can be represented as a linear combination of $\bm{b}_{1}, \dots, \bm{b}_{m-d}$.
So there exists a vector $\bm{y} \in \mathbb{R}^{m-d}$, such that $(\bm{b}_{1}, \dots, \bm{b}_{m-d})\bm{y}=\bm{b}_{m-d+1}$.

Let $\bm{s}$ be the vector of the coefficients of the input polynomials: $\bm{s}:=(f_{10}, \dots, f_{1m}, \dots, f_{n0}, \dots, f_{nm})$. 
Let the B\'ezout matrix of the input polynomials represented by $\bm{s}$ be $B(\bm{s})=(\bm{b}(\bm{s})_1, \dots, \bm{b}(\bm{s})_n)$.
In the same way, let the B\'ezout matrix of the polynomials we find ($\tilde{F}_1(x), \dots, \tilde{F}_n(x)$) be $B(\bm{s}+\Delta\bm{s})$, where $\Delta\bm{s}=(\tilde{f}_{10}-f_{10}, \dots, \tilde{f}_{1m}-f_{1m}, \dots, \tilde{f}_{n0}-f_{n0}, \dots, \tilde{f}_{nm}-f_{nm})$.

From the above, we can transfer Problem \ref{prob:AppGCD} to a constrained minimization problem with the objective function: 
\begin{equation}\label{eq:objectivefunction}
    \|\Delta\bm{s}\|=\sqrt{\sum_{i=1}^n\sum_{j=0}^m(\tilde{f}_{ij}-f_{ij})^2},
\end{equation}
and the constraints: 
\begin{equation}\label{eq:constraints}
    (\bm{b}(\bm{s}+\Delta\bm{s})_{1}, \dots, \bm{b}(\bm{s}+\Delta\bm{s})_{m-d})\bm{y}=\bm{b}(\bm{s}+\Delta\bm{s})_{m-d+1}\quad \text{for some} \quad \bm{y},
\end{equation}
with variables:
\begin{equation}\label{eq:variable}
    (\bm{\tilde{s}}, \bm{y})=(\tilde{f}_{10}, \dots, \tilde{f}_{1m}, \dots, \tilde{f}_{n0}, \dots, \tilde{f}_{nm}, y_1, \dots, y_{m-d}).
\end{equation}

\section{The Modified Newton Method}
\label{sec:newton}

We consider a constrained minimization problem of minimizing an objective function $f(\bm{x}):\mathbb{R}^s\rightarrow \mathbb{R}$
which is twice continuously differentiable, subject to the constraints 
$\bm{g}(\bm{x}) =(g_{1}(\bm{x}), \dots, g_{t}(\bm{x}))^T=\bm{0}$,
where $g_{i}(\bm{x})$ is a function of $\mathbb{R}^s\rightarrow \mathbb{R}$ and is also twice continuously differentiable.
For solving the constrained minimization problem, we use the modified Newton method by Tanabe (\cite{Tanabe1980}), which is a generalization of the Gradient Projection method (\cite{ros1961}),
as in the original GPGCD algorithm (\cite{TERUI2013127}).
For $\bm{x}_k$ which satisfies $\bm{g}(\bm{x}_k)=0$,
we calculate the search direction $\bm{d}_k$ and the Lagrange multipliers $\bm{\lambda}_k$ by solving the following linear system
\begin{equation}
\label{modifiedlinear}
\begin{pmatrix}
I & -(J_{\bm{g}}(\bm{x}_k))^T\\
J_{\bm{g}}(\bm{x}_k) & O
\end{pmatrix}
\begin{pmatrix}
\bm{d}_k\\
\bm{\lambda}_{k+1}
\end{pmatrix}
=-
\begin{pmatrix}
\nabla f(\bm{x}_k)\\
\bm{g}(\bm{x}_k)
\end{pmatrix},
\end{equation}
where $J_{\bm{g}}(\bm{x})$ is the Jacobian matrix represented as
\begin{equation}
J_{\bm{g}}(\bm{x})=\frac{\partial g_i}{\partial \bm{x}_j}.
\end{equation}

\section{The Algorithm for Calculating Approximate GCD}
\label{sec:algorithm}

We give an algorithm for calculating approximate GCD of multiple polynomials using the B\'{e}zout matrix in this section.

\subsection{Representation of the Jacobian matrix}

In the modified Newton method, the Jacobian matrix is represented as follows. 
From the constraints \eqref{eq:constraints} and the objective function \eqref{eq:objectivefunction}, 
let $d(i)=\deg(F_i(x)),i=1, \dots, n$, then the Jacobian matrix is represented as:

\begin{equation}
    \begin{pmatrix}
        \begin{array}{c|cccc|c}
        JL(2) & JM(2) & O & \dots & O & \\
        \vdots & O & \ddots & \ddots & \vdots & JR \\
        \vdots & \vdots & \ddots & \ddots & O & \\
        JL(n) & O & \dots & O & JM(n) & \\
        \end{array}
    \end{pmatrix},
    \label{eq:Jacobian}
\end{equation}
where
\begin{align*} 
    JL(k)&=JL_1(k)+JL_2(k), \\
    JL_1(k)&=(p_{1kij})_{i=1..m, j=1..m+1}, \\
    p_{1kij}&=
    \begin{cases}
        \sum_{l=0}^{\min(i, m-d-j+i)} f_{k(l)} y_{j-i+l} & 1 \le i < j \le m+1\\
        -\sum_{l=0}^{\min(m-i, m-j-d)}f_{k(i+l)} y_{j+l} & 1 \le j \le i \le m\\
    \end{cases},\\
    JL_2(k)&=(p_{2kij})_{i=1..m, j=1..m+1}, \\
    p_{2kij}&=
    \begin{cases}
        f_{k(m+1-d+i-j)} & 1 \le j \le i \le j+d-1 \le m\\
        -f_{k(m+1-d+i-j)} & 1 \le j-(m+1-d) \le i \le j-1 \le m\\
    \end{cases},\\
    JM(k)&=JM_1(k)+JM_2(k), \\
    JM_1(k)&=(p_{3kij})_{i=1..m, j=1..d(k)+1}, \\
    p_{3kij}&=
    \begin{cases}
        -\sum_{l=0}^{\min(i, \min(m-d,d(k))-j+i)} f_{1(l)} y_{j-i+l} & 1 \le i < j \le d(k)+1\\
        \sum_{l=0}^{\min(m-i, \min(m-d,d(k)-j)}f_{1(i+l)} y_{j+l} & 1 \le j \le i \le m\\
    \end{cases},\\
    JM_2(k)&=(p_{4kij})_{i=1..m, j=1..d(k)+1}, \\
    p_{4kij}&=
    \begin{cases}
        f_{1(m+1-d+i-j)} & 1 \le j \le i \le j+d-1 \le m\\
        -f_{1(m+1-d+i-j)} & 1 \le j-(m+1-d) \le i \le j-1 \le m\\
    \end{cases},\\
    JR&=(\bm{b}_{1}, \dots, \bm{b}_{m-d}).
\end{align*}

\subsection{Setting the Initial Values}

We give the initial values for variables in \eqref{eq:variable} as follows.
For the given polynomials \eqref{eq:defAppGCD}, we set the initial value for variables $\bm{s}_0$ as:
\begin{equation}
    \bm{s}_0=(f_{10}, \dots, f_{1m}, \dots, f_{n0}, \dots, f_{nm}).
\end{equation}
For the B\'ezout matrix $B=\Bez(F_1(x), \dots, F_n(x))=(\bm{b}_1, \dots, \bm{b}_m)$, we calculate $\bm{y}=(y_{01}, \dots, y_{0(m-d)})$ by solving the following least squares problem:
\begin{equation}
    \min\|(\bm{b}_1, \dots, \bm{b}_{m-d})\bm{y}-\bm{b}_{m-d+1}\|.
\end{equation}
Then we set the initial value for variables $\bm{y}_0$ as:
\begin{equation}\label{eq:initial_y}
    \bm{y}_0=(y_{01}, \dots, y_{0(m-d)}).
\end{equation} 
From above, we give the initial value for variables in \eqref{eq:variable} as
\begin{equation}
\label{eq:initialvalues}
(\bm{s}_0, \bm{y}_0)=(f_{10}, \dots, f_{1m}, \dots, f_{n0}, \dots, f_{nm}, y_{01}, \dots, y_{0(m-d)}).
\end{equation}

\subsection{Calculating the Approximate GCD and Finding the Approximate polynomials}

Let $(\bm{s}^*, \bm{y}^*) = (f_{10}^*, \dots, f_{1m}^*, \dots, f_{n0}^*, \dots, f_{nm}^*, y_{01}^*, \dots, y_{0(m-d)}^*)$ be the minimizer of the objective function in \eqref{eq:objectivefunction} calculated by the modified Newton method, corresponding to the coefficients of $F_1(x), \dots, F_n(x)$. 
Then, we calculate the approximate GCD from the B\'ezout matrix $B(\bm{s}^*)$ with \Cref{Theo:Barnett}.

Finally, we find the approximate polynomials $\tilde{F}_i(x)=\bar{F}_i(x)\tilde{H}(x)$, $i=1, \dots, n$ in \eqref{eq:Apppoly} by solving least squares problems:
\begin{equation}
    \label{eq:ls}
    \min\|F_i(x)-\bar{F}_i(x)\tilde{H}(x)\|\quad i=1, \dots, n.
\end{equation}

\subsection{The algorithm and running time analysis}
\label{sec:timeanalysis}

Summarizing above, we give the algorithm for calculating approximate GCD for multiple polynomials as Algorithm \ref{BezGPGCD}.
We give an analysis of the arithmetic running time of Algorithm \ref{BezGPGCD} as follows.

In Step \ref{alg:inivalue}, we set the initial values by the construction of the B\'ezout matrix, calculation of the initial values using the least squares solution, and the construction of the Jacobian matrix.
Since the dimension of the B\'ezout matrix and the Jacobian matrix is $mn \times m$ and $mn \times (m-d+\sum_{i=1}^n \deg(F_i)+1)$, respectively, 
we can estimate the running time of the construction of the B\'ezout matrix, calculation of the initial values using the least squares solution, 
and the construction of the Jacobian matrix as $O(m^2n)$ (\cite{Chionh2002}), $O(m^3)$ (\cite{doi:0.1007/978-0-387-40065-5}), $O(m(n(m-d+n)+\sum_{i=1}^n \deg(F_i)))=O(mn(m+n-d))$, respectively, where the Jacobian matrix is computed by \eqref{eq:Jacobian}.

In Step \ref{alg:linear}, since the dimension of the Jacobian matrix is $mn(m-d+\sum_{i=1}^n \deg(F_i)+1)$, the running time for solving the linear system is $O((mn+m-d+\sum_{i=1}^n \deg(F_i)+1)^3)=O((mn-d)^3)$ (\cite{doi:10.1137/1.9781611971446}).

In Step \ref{alg:condition}, the running time of the construction of the B\'ezout matrix and the Jacobian matrix is $O(m^2n)$ and $O(m(n(m-d+n)+\sum_{i=1}^n \deg(F_i)))=O(mn(m+n-d))$, respectively.

In Step \ref{alg:appgcd}, the running time for calculating the approximate GCD $\tilde{H}(x)$ and calculation of the approximate polynomials $\tilde{F}_i(x),~ i=1, \dots, n$ is $O(m^2n)$ and $O(n(m-d)^3)$.

As a consequence, the running time of Algorithm \ref{BezGPGCD} is the number of iteration times $O((mn-d)^3)$.

\begin{algorithm}\caption{The GPGCD algorithm for multiple polynomials with the B\'ezout matrix}
    \label{BezGPGCD}

    \noindent Inputs:
    \begin{itemize}
    \item[-]$F_1(x), \dots, F_n(x)\in\mathbb{R}(x)$: the given polynomials with $\max(\deg(F_i(x)), i=1\dots n) = \deg(F_1(x)) >0$,
    \item[-]$d \in \mathbb{N}$: the given degree of approximate GCD with $d \le \min(\deg(F_i(x)), i=1\dots n)$,
    \item[-]$\epsilon>0$: the stop criterion with the modified Newton method,
    \item[-]$0<\alpha\le1$: the step width with the modified Newton method.
    \end{itemize}
    Outputs:
    \begin{itemize}
    \item[-] $\tilde{H}(x)$: the approximate GCD, with $\deg(\tilde{H})=d$,
    \item[-] $\tilde{F}_1(x)$, \dots $\tilde{F}_n(x)$: the polynomials which are close to $F_1(x), \dots, F_n(x)$, respectively, with the GCD $\tilde{H}$.
    \end{itemize}
    \begin{enumerate}[Step 1]
    \item \label{alg:inivalue}Generate the B\'ezout matrix $B=\Bez(F_1(x) \dots F_n(x))$, the initial values $\bm{y}_0$ \eqref{eq:initial_y}, and the Jacobian matrix \eqref{eq:Jacobian}. 
    \item Set the initial values $(\bm{s}_0, \bm{y}_0)$ as in \eqref{eq:initialvalues}.
    \item \label{alg:linear}Solve the linear system \eqref{modifiedlinear} to find the search direction $\bm{d}_k$.\label{al:linsys}
    \item \label{alg:condition}If $\|\bm{d}_k\|<\epsilon$, obtain the $(\bm{s}^*, \bm{y}^*)$ as $(\bm{s}_k, \bm{y}_k)$, generate the B\'ezout matrix $\tilde{B}=\Bez(\tilde{F}_1(x), \dots, \tilde{F}_n(x))$ from polynomials $\tilde{F}_i(x),~ i=1, \dots, n$, then go to Step \ref{alg:appgcd}. Otherwise, let $(\bm{s}_{k+1}, \bm{y}_{k+1})=(\bm{s}_{k}, \bm{y}_{k})+\alpha\bm{d}_k$ and calculate the B\'ezout matrix and the Jacobian matrix with $(\bm{s}_{k+1}, \bm{y}_{k+1})$, then go to Step \ref{al:linsys}.
    \item \label{alg:appgcd} Calculate the approximate GCD $\tilde{H}(x)$ with \Cref{Theo:Barnett}. For $i=1,\dots,n$, calculate the approximate polynomials $\tilde{F}_i(x)$ by solving
    least squares problems in \eqref{eq:ls}. 
    Return $\tilde{F}_1(x),\dots,\tilde{F}_n(x)$ and $\tilde{H}(x)$.
    \end{enumerate}
\end{algorithm}

\section{Experiments}
\label{sec:experiments}

We have implemented our GPGCD algorithm on a computer algebra system 
Maple 2021.
For $i=1, \dots, n$, 
the test polynomials $F_i(x)$ 
are generated as 
\begin{equation}
\label{eq:experiments}
F_i(x)=\tilde{F}_i(x)H(x)+\frac{e}{\|\hat{F}_i(x)\|}\hat{F}_i(x). 
\end{equation}
Here, $\tilde{F}_i(x)$ 
and $H(x)$ are polynomials of degrees $\deg(F_i)-d$ 
and $d$, respectively.
Note that $F_i(x)$ 
are relatively prime polynomials.
They are generated as polynomials that their coefficients are floating-point numbers
and their absolute values are not greater than~10.\footnote{Coefficients are generated with the Mersenne Twister algorithm \cite{Matsumoto1998} by built-in function \texttt{Generate} with \texttt{RandomTools:-MersenneTwister} in Maple, which approximates a uniform distribution on $[-10, 10]$.}
The noise polynomials $\hat{F}_i(x)$ are polynomials of degree $\deg(F_i)-1$, which are randomly generated with coefficients given as the same as for $\tilde{F}_i(x)$ and $H(x)$.

We have generated 5 test groups of test polynomials, each group comprising 100 tests.
We set the degree of the input polynomials for all groups as 10, and the degree of the approximate GCD as $3, 4, 5, 6, 7$, respectively (shown as in Table \ref{table:tests}).
We set the number of polynomials for each tests as $n=10$, 
and set the norm of the noise $e$ as $e=0.01$ in our tests.
The stop criterion $\epsilon$ and the stop width $\alpha$ in Algorithm \ref{BezGPGCD} are set as $0.1$ and $1$, respectively.

We have carried out the tests on 
CPU Intel(R) Xeon(R) Silver 4210 at 2.20GHz with RAM 256GB under Linux 5.4.0 
with Maple 2021.

\subsection{The experimental results}

We have carried out the tests with the GPGCD algorithm with the B\'ezout matrix for multiple polynomials (abbreviated as GPGCD-B\'ezout-multiple algorithm).
For comparison, we have also carried out the tests with the original GPGCD algorithm with subresultant matrices for multiple polynomials (abbreviated as GPGCD-Sylvester-multiple algorithm) (\cite{Terui2010}).

For the results, we have compared the norm of the perturbation of the GPGCD-B\'ezout-multiple algorithm
with that of the GPGCD-Sylvester-multiple algorithm.
For $i=1, \dots, n$, let $F_{\textrm{B}i}=f_{\textrm{B}im}x^m+\dots+f_{\textrm{B}i0}x^0$ be the approximate polynomials from the outputs of the GPGCD-B\'ezout-multiple algorithm,
and $F_{\textrm{S}i}=f_{\textrm{S}im}x^m+\dots+f_{\textrm{S}i0}x^0$ be the approximate polynomials from the outputs of the GPGCD-Sylvester-multiple algorithm.
Then, let $\Delta B$ be the norm of the perturbation of the GPGCD-B\'ezout-multiple algorithm:
\begin{equation*}
    \Delta B=\sqrt{\sum_{i=1}^n\sum_{j=0}^m(f_{\textrm{B}ij}-f_{ij})^2},
\end{equation*}
and let $\Delta S$ be the norm of the perturbation of the GPGCD-Sylvester-multiple algorithm:
\begin{equation*}
    \Delta S=\sqrt{\sum_{i=1}^n\sum_{j=0}^m(f_{\textrm{S}ij}-f_{ij})^2}.
\end{equation*}
If we have 
    $|\Delta B - \Delta S| \le 0.5,$ 
then we say that the minimizers of both algorithms are \emph{closed to each other}.

The number of tests which the minimizers are closed to each other is shown in the left part of Table \ref{table:number}.
In `All', we show the number for all tests which the minimizers are closed to each other.
In `B\'ezout', we show the number for tests which the GPGCD-B\'ezout-multiple algorithm has a smaller perturbation than that of GPGCD-Sylvester-multiple algorithm.
In `Sylvester', we show the number for tests which GPGCD-Sylvester-multiple algorithm has a smaller perturbation than that of the GPGCD-B\'ezout-multiple algorithm.
In the same way, the number of tests which the minimizers are not closed to each other is shown in the right part of Table \ref{table:number}.

For the tests in which the minimizers are closed to each other, the average time and the average number of iterations are shown in Table \ref{table:time}.
In the left part, we show the average time for two algorithms.
In the middle part, we show the average number of iterations.
In the right part, we show the average time per iteration.

For the tests in which the minimizers are closed to each other, 
we have also calculated the norm of the remainders of the GPGCD-B\'ezout-multiple algorithm $\big{\|}(\|F_{\textrm{BR}1}(x)\|, \dots, \|F_{\textrm{BR}n}(x)\|)\big{\|}$, 
and the norm of the remainders of the GPGCD-Sylvester-multiple algorithm $\big{\|}(\|F_{\textrm{SR}1}(x)\|, \dots, \|F_{\textrm{SR}n}(x)\|)\big{\|}$, 
where
\begin{gather*}
    F_{\textrm{B}i}(x)=\bar{F}_i(x)\tilde{H}(x)+F_{\textrm{BR}i}(x), \quad \deg(F_{\textrm{BR}i}) < d, \quad i=1, \dots, n, \\
    F_{\textrm{S}i}(x)=\bar{F}_i(x)\tilde{H}(x)+F_{\textrm{SR}i}(x), \quad \deg(F_{\textrm{SR}i}) < d, \quad i=1, \dots, n. 
\end{gather*}
We show the average norm of the remainders in Table \ref{table:remainders}.

\subsection{Comparison with the GPGCD-Sylvester-multiple algorithm}

From Table \ref{table:number}, we see that the number of tests in which the minimizers are closed to each other of each group is larger when the degree of the approximate GCD is higher, where the least of them is over half of the total number of tests.
For all results, we see that the number of tests in which the perturbation of the GPGCD-Sylvester-multiple algorithm is smaller than that of the GPGCD-B\'ezout-multiple algorithm is larger than that of the opposite.
That means the GPGCD-Sylvester-multiple algorithm shows higher stability than the GPGCD-B\'ezout-multiple algorithm.

From Table \ref{table:time}, we see that total computing time, number of iterations, average computing time per iteration of the GPGCD-B\'ezout-multiple algorithm are smaller than those of the GPGCD-Sylvester-multiple algorithm, respectively.
By comparing average computing time per iteration in each group, we see that, except for Group 1 in the GPGCD-B\'ezout-multiple algorithm, when the degree of approximate GCD is higher, the average time per iteration of both algorithms is smaller. 

In the previous research, the GPGCD algorithm with the B\'ezout matrix with 2 polynomials (abbreviated as GPGCD-B\'ezout algorithm) (\cite{10.1007/978-3-030-60026-6_10}), we have shown that the average of the norm of the remainders of the tests in the GPGCD-B\'ezout algorithm is larger than that of the original GPGCD algorithm (abbreviated as GPGCD-Sylvester algorithm).
From Table \ref{table:remainders}, we see that there are little difference between the GPGCD-B\'ezout-multiple algorithm and the GPGCD-Sylvester-multiple algorithm in the sense of the average of the norm of the remainders.


\begin{table}
    \centering
    \caption{Degrees of test polynomials \eqref{eq:experiments}}
    \begin{tabular}{ccccccccc}\hline
    \label{table:tests}
    Group & $m=\deg(F_i)$ & $d=\deg(H)$ \\ \hline
    1 & 10 & 3 \\ \hline
    2 & 10 & 4 \\ \hline
    3 & 10 & 5 \\ \hline
    4 & 10 & 6 \\ \hline
    5 & 10 & 7 \\ \hline

    \end{tabular}
\end{table}

\begin{table}
    \centering
    \caption{The number of tests which the minimizers are closed to / not closed to each other}
    \begin{tabular}{ccccccccc}\hline
    \label{table:number}
     & \multicolumn{6}{c}{The number of tests}\\ \cline{2-7}
    Group & \multicolumn{3}{c}{Closed to each other} & \multicolumn{3}{c}{Not closed to each other} \\ \cline{2-7}
     & All & B\'ezout & Sylvester & All & B\'ezout & Sylvester \\ \hline
    1 & 55 & 28 & 27 & 45 & 4 & 41 \\ \hline
    2 & 57 & 21 & 36 & 43 & 9 & 34 \\ \hline
    3 & 72 & 29 & 43 & 28 & 2 & 26 \\ \hline
    4 & 72 & 24 & 48 & 28 & 0 & 28 \\ \hline
    5 & 81 & 19 & 62 & 19 & 2 & 17 \\ \hline
    \end{tabular}
\end{table}

\begin{table}
    \centering
    \caption{Comparison of computing time and the number of iterations (See Remark 1 for detail)}
    \begin{tabular}{ccccccccccccc}\hline
    \label{table:time}
     & \multicolumn{2}{c}{Average time (sec.)} & \multicolumn{2}{c}{Average number of iterations} & \multicolumn{2}{c}{Average time per iteration (sec.)} \\ \cline{2-7}
    Group & B\'ezout & Sylvester & B\'ezout & Sylvester & B\'ezout & Sylvester \\ \hline
    1 & 9.804 & 166.266 & 14.47 & 23.47 & 0.677 & 7.083 \\ \hline
    2 & 11.035 & 141.235 & 16.14 & 21.28 & 0.684 & 6.637 \\ \hline
    3 & 9.587 & 112.952 & 14.375 & 18.14 & 0.667 & 6.227 \\ \hline
    4 & 7.296 & 99.230 & 11.33 & 17 & 0.644 & 5.837 \\ \hline
    5 & 6.272 & 89.058 & 9.83 & 16.09 & 0.638 & 5.536 \\ \hline
    \end{tabular}
\end{table}

\begin{table}
    \centering
    \caption{Comparison of norm of the remainders}
    \begin{tabular}{ccccccccccc}\hline
    \label{table:remainders}
    Group & B\'ezout & Sylvester \\ \hline
    1 & $2.320 \times 10^{-9}$ & $1.875 \times 10^{-9}$ \\ \hline
    2 & $1.390 \times 10^{-9}$ & $9.170 \times 10^{-9}$ \\ \hline
    3 & $1.150 \times 10^{-9}$ & $1.059 \times 10^{-9}$ \\ \hline
    4 & $9.684 \times 10^{-9}$ & $8.427 \times 10^{-10}$ \\ \hline
    5 & $7.987 \times 10^{-11}$ & $9.012 \times 10^{-11}$ \\ \hline
    \end{tabular}
\end{table}

\section{Conclusions}
We have proposed an algorithm using the B\'ezout matrix based on the GPGCD algorithm for multiple polynomials.

The experimental results show that, in most of the cases, the minimizer of the proposed algorithm is closed to the minimizer of the GPGCD-Sylvester-multiple algorithm. 
For the cases in which the minimizers of both algorithms are closed to each other, 
total computing time, the number of iterations, average computing time per iteration of the GPGCD-B\'ezout-multiple algorithm are smaller than those of the GPGCD-Sylvester-multiple algorithm, respectively.
Thus, we see that the proposed algorithm is more efficient than the GPGCD-Sylvester-multiple algorithm.
By comparing the computing time per iteration, the experimental results show that when the degree of the approximate GCD is higher, the computing time per iteration of the proposed algorithm is smaller. This is consistent with the arithmetic running time analysis shown in Section \ref{sec:timeanalysis}.

The previous research has shown that the GPGCD-B\'ezout algorithm was not so accurate in the sense of the norm of perturbations compared with the GPGCD-Sylvester algorithm. 
However, new experimental results show that, with the improvement for calculating the approximate polynomials as shown in \eqref{eq:ls}, 
the GPGCD-B\'ezout-multiple algorithm calculate approximate GCDs with almost the same accuracy
as the GPGCD-Sylvester-multiple algorithm in the sense of the norm of perturbations, for the tests in which the minimizers of both algorithms are closed to each other.

On the other hand, the experimental results show that the number of tests in which the perturbation of the GPGCD-Sylvester-multiple algorithm is smaller than that of the proposed algorithm is larger than that of the opposite.
Thus, the GPGCD-Sylvester-multiple algorithm shows higher stability than the proposed algorithm. 
Improving stability of the proposed algorithm will be one of the topics of further research.

\section*{Acknowledgments}

This research was partially supported by JSPS KAKENHI Grant Number JP20K11845.
The author (Chi) was supported by the research assistant fellowship from Graduate School of 
Pure and Applied Sciences, University of Tsukuba.



\end{document}